\documentclass{svproc}
\usepackage{url}

\def\be{\begin{equation}}
\def\en{\end{equation}}
\def\bee{\begin{eqnarray*}}
\def\ene{\end{eqnarray*}}
\def\R{{\bf R}}
\def\P{{\bf P}}
\def\E{{\bf E}}
\def\H{{\bf H}}
\def\Var{{\rm Var}}

\def\ep{\varepsilon}
\begin{document}
\mainmatter              % start of a contribution
\title{Two-sided inequalities \\ for the density function's maximum\\ of weighted sum of chi-square variables}
\titlerunning{Two--sided inequalities for PDF's maximum}  % abbreviated title (for running head)
\author{Sergey G. Bobkov\inst{1, 2}  \and Alexey A. Naumov\inst{2} \and \\
Vladimir V. Ulyanov\inst{2, 3} }
\authorrunning{S. Bobkov, A. Naumov,
V. Ulyanov} % abbreviated author list (for running head)
\tocauthor{Sergey Bobkov, Alexey Naumov,
Vladimir Ulyanov}
\institute{University of Minnesota, \\  Vincent Hall 228, 206 Church St SE, Minneapolis, MN 55455 USA,\\
%\\ WWW home page:
%\texttt{http://users/\homedir iekeland/web/welcome.html}
\and
Faculty of Computer Science, \\HSE University,  109028 Moscow, Russian Federation,
\and 
Faculty of Computational Mathematics and Cybernetics, \\Lomonosov Moscow State
University,\\ 119991 Moscow, Russian Federation\\
\email{vulyanov@cs.msu.ru}
}

\maketitle              % typeset the title of the contribution

\begin{abstract}
Two--sided bounds are constructed for a probability density function of a weighted sum of  chi-square variables. Both cases of central and non-central chi-square variables are considered. The upper and lower bounds have the same dependence on the parameters of the sum and differ only in absolute constants. The estimates obtained will be useful, in particular, when comparing two Gaussian random elements in a Hilbert space and in multidimensional central limit theorems, including the infinite-dimensional case.  
% We would like to encourage you to list your keywords within
% the abstract section using the \keywords{...} command.
\end{abstract}
\section{Introduction}
In many statistical and probabilistic applications, we have to solve the
problem of Gaussian comparison, that is, one has
to evaluate how the probability of a ball under a Gaussian measure is affected,
if the mean and the covariance operators of this Gaussian measure are
slightly changed. In \cite{GNSU19} 
 we present particular examples motivating the results when such
``large ball probability'' problem naturally arises,
including bootstrap validation, Bayesian inference and high-dimensional CLT, see also \cite{PU13} and  \cite{FU20}.
%This paper presents sharp bounds for the Kolmogorov distance between
%the probabilities of two Gaussian elements to hit a ball in a Hilbert space.
%The key property of these bounds is that they are dimension-free and
%depend on the nuclear (Schatten-one) norm of the difference between the
%covariance operators of the elements.
% 
%We also state a tight dimension free anti-concentration bound for a
%squared norm of a Gaussian element in Hilbert space which refines the
%well-known results on the density
%of a chi-squared distribution; see Theorem~\ref{band of GE}.
The tight non-asymptotic bounds for the Kolmogorov distance
between the probabilities of two Gaussian elements to hit a ball in a
Hilbert space have been derived in \cite{GNSU19} and  \cite{NSTU18}.
The key property of these bounds is that they are dimension-free and
depend on the nuclear (Schatten-one) norm of the difference
between the covariance operators of the elements and on the norm of the
mean shift.
The obtained bounds significantly improve the bound based on Pinsker's
inequality via the Kullback--Leibler divergence.
It was also established an anti-concentration bound for a squared norm $||Z - a||^2, \ \ a \in \H,$ of a shifted 
 Gaussian element $Z$ with zero mean in a Hilbert space $\H$. The decisive role in proving the results was played by the upper estimates for the maximum of the probability density function $g(x, a)$ of $||Z - a||^2$, see Theorem 2.6 in \cite{GNSU19}:
\be
\label{Be}
\sup_{x\geq 0} g(x, a) \leq  c\, (\Lambda_1 \Lambda_2)^{-1/4},
\en
where $c$ is an absolute constant and
$$
\Lambda_1 = \sum_{k=1}^{\infty}  \lambda_k^2, \qquad  
\Lambda_2 = \sum_{k=2}^{\infty}  \lambda_k^2
$$ 
with $ \lambda_1 \geq   \lambda_2 \geq \dots$ are the eigenvalues of a covariance operator $\Sigma$ of $Z$.

It is well known that $ g(x, a)$ can be considered as a density function of a weighted sum
of non-central $ \chi^{2}$ distributions.
An explicit but cumbersome representation for $ g(x, a) $ in
finite dimensional space $\H$ is available (see, e.g., Section~18 in
Johnson, Kotz and Balakrishnan \cite{JKotzB1994}).
However, it involves some special characteristics of the related
Gaussian measure which makes it hard to use
in specific situations.
Our result  (\ref{Be})
is much more transparent and provide
sharp uniform upper bounds.
Indeed, in the case $\H= \R^{d}$, $a = 0$,  $\Sigma$ is the unit matrix, one has that
the distribution of $||Z||^{2}$ is the standard $\chi^{2}$ with
$d$ degrees of freedom and the maximum of its probability density function is proportional to $d^{-1/2}$. This is the same as what we get in (\ref{Be}).
%%
%\begin{eqnarray*}
%\label{identity case} \sup_{x \geq0} p_{\xiv}(x, 0) \asymp
%p^{-1/2}.
%\end{eqnarray*}
%%
%Hence, the bound \eqref{density bound 0} gives the right dependence on
%$\dimp$ because $\CONSTdlt(\Sigma_{\xiv}) \asymp p^{-1/2}$.
%%See also Example after Theorem~\ref{band of GE}.
%However, a lower bound for $\sup_{x \geq0} p_{\xiv}(x, \av) $ in
%the general case is still an open question.
%One can even get two-sided bounds for $p_{\xiv}(x, \av)$ but under
%additional conditions, see, for example, Christoph, Prokhorov and
%Ulyanov \cite{christoph1996}.

At the same time,  
 it was noted in \cite{GNSU19} that obtaining lower estimates for $\sup_x g(x, a)$ remains an open problem. 
The latter problem was partially solved in \cite{christoph1996}, Theorem~1. However, it was done under additional conditions and we took into account the multiplicity of the largest
eigenvalue.

In the present paper we get two--sided bounds for $\sup_x g(x, 0)$ in the finite-dimensional case $\H = \R^d$, see Theorem~\ref{1} below. The bounds are dimension-free, that is they do not depend on $d$.  Thus, for the upper bounds (\ref{Be}), we obtain a new proof, which is of independent interest.  And new lower bounds show the optimality of (\ref{Be}),
%the upper bounds in \cite{GNSU19}, 
since the upper and lower bounds differ only in absolute constants. Moreover, new two-sided bounds are constructed for $\sup_x g(x, a)$  with  $a\neq 0$ in the finite-dimensional case $\H = \R^d$, see Theorem~\ref{2} below. Here we consider a typical situation, where $\lambda_1$ does not  
  dominate the other coefficients.

\section{Main results}
\noindent
For independ standard normal random variables $Z_k \sim N(0,1)$, consider the weighted sum
$$
W_0 = \lambda_1 Z_1^2 + \dots + \lambda_n Z_n^2, \qquad 
\lambda_1 \geq \dots \geq \lambda_n > 0.
$$
It has a continuous probability density function $p(x)$ on the positive half-axis. 
Define the functional
$$
M(W_0) = \sup_x \, p(x).
$$

%\vskip2mm
\begin{theorem}
\label{1}
Up to some absolute constants $c_0$ and $c_1$, we have
\be
\label{th1}
c_0 (A_1 A_2)^{-1/4} \leq M(W_0) \leq c_1 (A_1 A_2)^{-1/4},
\en
where
$$
A_1 = \sum_{k=1}^n  \lambda_k^2, \qquad  
A_2 = \sum_{k=2}^n  \lambda_k^2
$$
and
$$
c_0 = \frac{1}{4e^2\sqrt{2\pi}} > 0.013, \qquad  
c_1 = \frac{2}{\sqrt{\pi}} < 1.129.
$$
\end{theorem}

Theorem 1 can be extended to more general weighted sums:
\be
\label{non-zero}
W_a = \lambda_1 (Z_1 - a_1)^2 + \dots + \lambda_n (Z_n - a_n)^2 \nonumber
\en
with parameters $\lambda_1 \geq \dots \geq \lambda_n > 0$ and
$a = (a_1, \dots, a_n) \in \R^n$. 

It has a continuous probability density function $p(x, a)$ on the positive half-axis $x > 0$.  
Define the functional
$$
M(W_a) = \sup_x \, p(x, a).
$$

{\bf Remark.} It is known that for any non-centred Gaussian element $Y$ in a Hilbert space, the random variable $||Y||^2$ is distributed as $\sum_{i=1}^{\infty}\lambda_i (Z_i - a_i)^2$ with some real $a_i$ and $\lambda_i$ such that 
$$\lambda_1 \geq \lambda_2 \geq \dots \geq 0 \ \ \ {\rm and}  \ \ \ \sum_{i=1}^{\infty}\lambda_i < \infty .$$ 
Therefore, the upper bounds for $M(W_a)$ immediately imply the upper bounds for the probability density function of $||Y||^2$. 

\smallskip

\begin{theorem}
\label{2}
If  $\lambda_1^2 \leq %\frac{1}{3}
 A_1/3$, then %for the maximum of the random variable $W$ defined in (6), 
one has a two-sided bounds
$$
\frac{1}{4\sqrt{3}}\,\frac{1}{\sqrt{A_1 + B_1}} \leq
M(W_a) \leq \frac{2}{\sqrt{A_1 + B_1}},
$$
where
$$
A_1 = \sum_{k=1}^n  \lambda_k^2, \qquad  
B_1 = \sum_{k=1}^n  \lambda_k^2 a_k^2.
$$
Moreover, the left inequality holds  without any assumptions on $\lambda_1^2$.
\end{theorem}

{\bf Remark.} In Theorem \ref{2}  we only consider a typical situation, where $\lambda_1$ does not  
  dominate the other coefficients. Moreover,  
 the condition $\lambda_1^2 \leq 
%\frac{1}{3} 
A_1/3$ necessarily
implies that $n \geq 3$.
If this condition is violated, 
 %(that is, for $\lambda_1^2 > \frac{1}{3} A_1$),
the behaviour of $M(W_a)$ should be studied separately.

\section{Auxiliary results}

For the lower bounds in the theorems, one may apply the following lemma, which goes back to 
the work by Statulyavichus \cite{St65}, see also Proposition 2.1 in \cite{BC15}.
%general elementary inequality
\begin{lemma}\label{low}
Let  $\eta$ be a random variable with $M(\eta)$ denoting the maximum of its probability density function. Then one has
\be
\label{lower}
M^2(\eta) \,\Var(\eta) \geq \frac{1}{12}.
\en
Moreover,  the equality in (\ref{lower}) is attained 
for the uniform distribution on any finite interval.
\end{lemma}
{\bf Remark.} There are multidimensional extensions of (\ref{lower}), see e.g. \cite{B88}, \cite{H80} and Section III in \cite{B11}.
%of this inequality.

\smallskip
\noindent
{\bf Proof.} Without loss of generality we may assume that $M(\eta) =1.$ 

Put $H(x) = \P (|\eta - \E \eta| \geq x), \quad x\geq 0.$ 

Then, $H(0) = 1$ and $H'(x)\geq -2$, which gives $H(x)\geq 1 - 2x,$ so
\begin{eqnarray}
\Var(\eta) &=& 2\int_0^{\infty} xH(x)\, dx \geq 2\int_0^{1/2} xH(x)\, dx \nonumber\\
&\geq & 2\int_0^{1/2} x(1-2x)\, dx =\frac{1}{12}.\nonumber
\end{eqnarray} 
Lemma is proved.

\vskip3mm 

The following lemma will give the lower bound in Theorem~\ref{2}.

\begin{lemma}\label{lowg} For the random variable $W_a$ defined in (\ref{non-zero}), the maximum $M(W_a)$ of its probability density function satisfies
%\be
%M(W) \geq \frac{1}{4\sqrt{2}}\,\frac{1}{\sqrt{A_1 + B_1}}.
%\en
\be
\label{lowerg}
M(W_a) \geq \frac{1}{4\sqrt{3}}\,\frac{1}{\sqrt{A_1 + B_1}},
\en
where 
$$
A_1 = \sum_{k=1}^n  \lambda_k^2, \qquad  
B_1 = \sum_{k=1}^n  \lambda_k^2 a_k^2.
$$

\end{lemma}
{\bf Proof.} 
Given $Z \sim N(0,1)$ and $b \in \R$, we have
$$
\E\,(Z - b)^2 = 1 + b^2, \qquad 
\E\,(Z - b)^4 = 3 + 6b^2 + b^4,
$$
so that $\Var((Z - b)^2) = 2 + 4b^2$. It follows that
$$
\Var(W_a) = \sum_{k=1}^n \lambda_k^2\,(2 + 4a_k^2) = 
2 A_1 + 4B_1 \leq 4 (A_1 + B_1).
$$

Applying (\ref{lower}) with $\eta = W_a$, we arrive at (\ref{lowerg}).

Lemma is proved.

\smallskip

The proofs of the upper bounds in the theorems are based on the following lemma.

\begin{lemma}\label{L1}
Let 
$$\alpha_1^2 + \dots + \alpha_n^2 = 1.$$ 
If $\alpha_k^2 \leq {1}/{m}$ for $m = 1,2,\dots$, 
then the characteristic function $f(t)$ of  
the random variable
$$
W = \alpha_1 Z_1^2 + \dots + \alpha_n Z_n^2
$$
satisfies
\be
\label{char}
|f(t)| \leq \frac{1}{(1 + 4t^2/m)^{m/4}}.
\en
In particular, in the cases $m=4$ and $m = 3$, $W$ has a bounded 
density with $M(W) \leq {1}/{2}$ and $M(W) < 0.723$ respectively.
\end{lemma}
{\bf Proof.} Necessarily $n \geq m$. The characteristic function has the form
$$
f(t) = \prod_{k=1}^n (1 - 2\alpha_k it)^{-1/2},
$$
so
$$
- \log |f(t)| = \frac{1}{4} \sum_{k=1}^n \log(1 + 4\alpha_k^2 t^2).
$$

First, let us describe the argument in the simplest case $m=1$.

For a fixed $t$, consider the concave function
$$
V(b_1,\dots,b_n) = \sum_{k=1}^n \log(1 + 4b_k t^2)
$$
on the simplex
$$
Q_1 
 \, = \,
\Big\{(b_1,\dots,b_n):  b_k \geq 0, \ b_1 + \dots + b_n = 1\Big\}.
$$
It has $n$ extreme points $b^k = (0,\dots,0,1,0,\dots,0)$. Hence
$$
\min_{b \in Q_1} V(b) = V(b^{k}) = \log(1 + 4t^2),
$$
that is, $|f(t)| \leq (1 + 4t^2)^{-1/4}$, which corresponds to (\ref{char}) for $m=1$.

If $m=2$, we consider the same function $V$ on the convex set
$$
Q_2 = \Big\{(b_1,\dots,b_n): 
0 \leq b_k \leq \frac{1}{2}, \ b_1 + \dots + b_n = 1\Big\},
$$
which is just the intersection of the cube $[0,\frac{1}{2}]^n$ with the 
hyperplane. It has ${n(n-1)}/{2}$ extreme points $$b^{kj}, \ \ \   
1 \leq k < j \leq n,$$ with coordinates ${1}/{2}$ on the $j$-th and $k$-th 
places and with zero elsewhere. Indeed, suppose that a point 
$$b = (b_1,\dots,b_n) \in Q_2$$ 
has at least two non-zero coordinates 
$0 < b_k, b_j < {1}/{2}$ for some $k < j$. Let $x$ be the point 
with coordinates 
$$x_l = b_l \ \ \ {\rm for}  \ \ \  l \neq k,j,\,\,\, x_k = b_k + \ep, \ \ \ 
{\rm and} \ \ \  
 x_j = b_j - \ep, $$
and similarly, let $y$ be the point such that 
$$ y_l = b_l \ \ \  {\rm for} \ \ \  l \neq k,j, \,\,\, y_k = b_k - \ep, \ \ \   {\rm and} \ \ \  y_j = b_j + \ep.$$
If $\ep>0$ is small enough, then both $x$ and $y$ lie in $Q_2$, while 
$$b = (x+y)/2, \ \ \  x \neq y.$$ 
Hence such $b$ cannot be an extreme point. 
Equivalently, any extreme point $b$ of $Q_2$ is of the form $$b^{kj}, \ \ \   
1 \leq k < j \leq n.$$ 
Therefore, we  conclude that
$$
\min_{b \in Q_2} V(b) = V(b^{kj}) = 2\log(1 + 2t^2),
$$
which is the first desired claim. 

In the general case, consider the function $V$
on the convex set
$$
Q_m = \Big\{(b_1,\dots,b_n): 
0 \leq b_k \leq \frac{1}{m}, \ b_1 + \dots + b_n = 1\Big\}.
$$
By a similar argument, any extreme point $b$ of $Q_m$ has zero 
for all coordinates except for $m$ places where the coordinates are 
equal to ${1}/{m}$. Therefore,
$$
\min_{b \in Q_m} V(b) = 
V\Big(\frac{1}{m},\dots,\frac{1}{m},0,\dots,0\Big) = 
m\log(1 + 4t^2/m),
$$
and we are done. 

In case $m=4$, using the inversion formula, we get
$$
M(W) \, \leq \, 
\frac{1}{2\pi} \int_{-\infty}^\infty |f(t)|\,dt 
 \, \leq \,
\frac{1}{2\pi} \int_{-\infty}^\infty \frac{1}{1 + t^2}\,dt
 \, = \,
\frac{1}{2}.
$$
Similarly, in the case $m=3$,
$$
M(W) \, \leq \, \frac{1}{2\pi} \int_{-\infty}^\infty 
\frac{1}{(1 + \frac{4}{3}\,t^2)^{3/4}}\,dt
 \, < \, 0.723.
$$
Lemma is proved.

\section{Proofs of main results}

{\bf Proof of Theorem \ref{1}.} 
%\vskip5mm
In the following we shall write $W$ instead of $W_0$.

If $n=1$,  then 
the distribution function and the probability density function of $W = \lambda_1 Z_1^2$ are given by
$$
F(x) = 2\,\Phi\bigg(\sqrt{\frac{x}{\lambda_1}}\,\bigg) - 1, \quad
p(x) = \frac{1}{\sqrt{2\pi \lambda_1}}\,e^{-x/(2\lambda_1)} \qquad (x > 0),
$$
respectively. 
Therefore, $p$ is unbounded near zero, so that $M(W) = \infty$. This is 
consistent with (\ref{th1}), in which case $A_1 = \lambda_1^2$ and $A_2 = 0$.

If $n=2$, the density $p(x)$ is described as the convolution
\be
\label{conv}
p(x) \, = \, \frac{1}{2\pi \sqrt{\lambda_1 \lambda_2}} \int_0^1
\frac{1}{\sqrt{(1-t) t}}\,\exp\Big\{-\frac{x}{2}\,
\Big[\frac{1-t}{\lambda_1} + \frac{t}{\lambda_2}\Big]\Big\}\,dt \qquad (x>0). 
\en
Hence, $p$ is decreasing and attains maximum at $x=0$:
$$
M(W) = \frac{1}{2\pi \sqrt{\lambda_1 \lambda_2}} \int_0^1
\frac{1}{\sqrt{(1-t) t}}\,dt = \frac{1}{2\sqrt{\lambda_1 \lambda_2}}.
$$
Since $A_1 = \lambda_1^2 + \lambda_2^2$ and $A_2 = \lambda_2^2$, 
we conclude, using the assumption $\lambda_1 \geq \lambda_2$, that
$$
\frac{1}{2}\,(A_1 A_2)^{-1/4} \leq
M(W) \leq \frac{1}{2^{3/4}}\,(A_1 A_2)^{-1/4}.
$$

As for the case $n \geq 3$, the density $p$ is vanishing at zero and
attains maximum at some point $x>0$. 
%It would be interesting to verify whether or not $p$ is unimodal.

The furher proof of Theorem \ref{1} is based on the following observations and Lemma~\ref{L1}.

By homogeneity of (\ref{th1}), we may assume 
that $A_1 = 1$. %and also that $n \geq 5$. 

If $\lambda_1 \leq {1}/{2}$, 
then all $\lambda_k^2 \leq {1}/{4}$, so that $M(W) \leq {1}/{2}$, 
by Lemma \ref{L1}. Hence, the inequality of the form
$$
M(W) \leq \frac{1}{2}\,(A_1 A_2)^{-1/4} 
$$
holds true. %with $c = \frac{1}{2}$ like in case $n=2$.

Now, let $\lambda_1 \geq {1}/{2}$, so that $A_2 \leq {3}/{4}$. Write
$$
W = \lambda_1 Z_1^2 + \sqrt{A_2}\,\xi, \qquad \xi = 
\sum_{k=2}^n \alpha_k Z_k^2, 
\quad \alpha_k = \frac{\lambda_k}{\sqrt{A_2}}.
$$
By construction, $\alpha_2^2 + \dots + \alpha_n^2 = 1$.

{\it Case} 1: $\lambda_2 \geq %\frac{1}{2}
\sqrt{A_2}/2$. Since the function 
$M(W)$ may only decrease when adding an independent random variable 
to $W$, we get using (\ref{conv}) that
$$
M(W) \leq M(\lambda_1 Z_1^2 + \lambda_2 Z_2^2)
= \frac{1}{2\sqrt{\lambda_1 \lambda_2}} \leq c\,(A_1 A_2)^{-1/4},
$$
where the last inequality holds with $c = 1$.
This gives the upper bound in (\ref{th1}) with constant 1.

{\it Case} 2: $\lambda_2 \leq %\frac{1}{2}
\sqrt{A_2}/2$. It implies that $n\geq 5$ and all 
$\alpha_k^2 \leq {1}/{4}$ for $k > 1$. By Lemma~\ref{L1} with $m=4$, the random variable $\xi$ 
has the probability density function $q$ bounded by ${1}/{2}$. The distribution function
of %the random variable 
$W$ may be written as
$$
\P\{W \leq x\} = 
\int_0^{x/\sqrt{A_2}} \P\Big\{|Z_1| \leq \frac{1}{\sqrt{\lambda_1}}\,
(x - y\sqrt{A_2})^{1/2}\Big\}\,q(y)\,dy, \quad x > 0,
$$
and its density has the form
$$
p(x) = \frac{1}{\sqrt{2\pi \lambda_1}} \int_0^{x/\sqrt{A_2}} 
\frac{1}{\sqrt{x - y\sqrt{A_2}}}\,e^{-(x - y\sqrt{A_2})/(2\lambda_1)}\,q(y)\,dy.
$$
Equivalently,
\be
\label{convlow}
p(x\sqrt{A_2}) = \frac{1}{\sqrt{2\pi \lambda_1}}\, A_2^{-1/4}
\int_0^{x} \frac{1}{\sqrt{x - y}}\,e^{-(x - y)\sqrt{A_2}/(2\lambda_1)}\,q(y)\,dy.
\en
Since $\lambda_1 \geq {1}/{2}$, we immediately obtain that
$$
M(W) \leq A_2^{-1/4}\,\frac{1}{\sqrt{\pi}}\ \sup_{x>0} \ 
\int_0^x \frac{1}{\sqrt{x - y}}\,q(y)dy.
$$
But, using $q \leq {1}/{2}$, we get
\bee
\int_0^x \frac{1}{\sqrt{x - y}}\,q(y)dy 
 & = & 
\int_{0 < y < x, \ x-y < 1} \frac{1}{\sqrt{x - y}}\,q(y)dy \\
&&+
\int_{0 < y < x, \ x-y > 1} \frac{1}{\sqrt{x - y}}\,q(y)dy \\
 & \leq &
\frac{1}{2}\,\int_0^1 \frac{1}{\sqrt{z}}\,dz + 1 \, = \, 2.
\ene
Thus,
$$
M(W) \leq 2 A_2^{-1/4}\,\frac{1}{\sqrt{\pi}}.
$$
Combining the obtained upper bounds for $M(W)$ in all cases we get 
 the upper bound in (\ref{th1}).

For the lower bound, one may apply the inequality (\ref{lower}) in Lemma~\ref{low}.
%\be
%M(\eta)^2\,\Var(\eta) \geq \frac{1}{12},
%\en
%where $\eta$ is an arbitrary random variable with density. It goes back to 
%the work by Statulyavichus (1965); the equality in this relation is attained 
%for the uniform distribution on any finite interval. 
Thus, we obtain that
$$
M(W) \geq \frac{1}{2\sqrt{6}}
$$
due to the assumption $A_1 = 1$ and the property $\Var(Z_1^2) = 2$.

If $\lambda_1^2 \leq {1}/{2}$, we have $A_2 \geq {1}/{2}$.  Hence,
\be
\label{th1low1}
M(W) \geq \frac{1}{2\sqrt{6}} \geq c_0\,(A_1 A_2)^{-1/4},
\en
where the last inequality holds true with 
$$c_0 = \frac{1}{2^{5/4}\sqrt{6}} \, \geq 0.171.$$
%(\frac{1}{2})^{1/4} = 0.121...$$

In case $\lambda_1^2 \geq \frac{1}{2}$, we have $A_2 \leq {1}/{2}$. 
Returning to the formula (\ref{convlow}),
let us choose $x = \E \xi + 2$ and restrict the integration to the interval
$$
\Delta: \max(\E \xi - 2,0) < y < \E \xi + 2.
$$
On this interval necessarily $$x - y \leq 4.$$ 
Therefore, (\ref{convlow}) yields
$$
M(W) \geq \frac{A_2^{-1/4}}{2\sqrt{2\pi \lambda_1}} \,\cdot
e^{-2\sqrt{A_2}/\lambda_1}\,\P\{\xi \in \Delta\}.
$$
Here, 
$$
\frac{A_2}{\lambda_1^2} = \frac{1}{\lambda_1^2} - 1 \leq 1,
$$
and we get
$$
M(W) \geq \frac{A_2^{-1/4}}{2\sqrt{2\pi}} \,\cdot e^{-2}\,\P\{\xi \in \Delta\}.
$$

Now, recall that $\xi \geq 0$ and $\Var(\xi) = 
2\,(\alpha_2^2 + \dots + \alpha_n^2) = 2$. 
Hence, by Chebyshev's inequality,
$$
\P\{|\xi - \E\xi| \geq 2\} \leq 
\frac{1}{4}\,\Var(\xi) = \frac{1}{2}.
$$
That is, $\P\{\xi \in \Delta\} \geq {1}/{2}$, and thus
$$
M(W) \geq \frac{ (A_1 A_2)^{-1/4} }{4\sqrt{2\pi}} \,e^{-2} \geq 0.013\cdot  (A_1 A_2)^{-1/4} .
$$
Theorem \ref{1} is proved.

{\bf Proof of Theorem \ref{2}.}  In the following we shall write $W$ instead of $W_a$.

The lower bound in Theorem \ref{2} immediately follows from (\ref{lowerg}) in Lemma \ref{lowg} without any assumption on $\lambda_1^2.$
%We only consider a typical situation where none of the coefficients
%$\lambda_k$ dominates the other coefficients. In order to make 
%a reasonable hypothesis about the magnitude of the maximum 
%of the density of $W$, let us return to the general lower 
%bound (5). Given $Z \sim N(0,1)$ and $a \in \R$, we have
%$$
%\E\,(Z - a)^2 = 1 + a^2, \qquad 
%\E\,(Z - a)^4 = 3 + 6a^2 + a^4,
%$$
%so that $\Var((Z - a)^2) = 2 + 4a^2$. It follows that
%$$
%\Var(W) = \sum_{k=1}^n \lambda_k^2\,(2 + 4a_k^2) = 
%2 A_1 + 4B_1,
%$$
%where
%$$
%A_1 = \sum_{k=1}^n  \lambda_k^2, \qquad  
%B_1 = \sum_{k=1}^n  \lambda_k^2 a_k^2.
%$$
%Applying (5) with $\eta = W$, we arrive at the following lower bound.
%
%\vskip5mm
%{\bf Lemma 3.} For the random variable $W$ as in (6), 
%\be
%M(W) \geq \frac{1}{4\sqrt{2}}\,\frac{1}{\sqrt{A_1 + B_1}}.
%\en

%\vskip5mm
Our next aim is to reverse this bound up to a numerical factor
under suitable natural assumptions. 

Without loss of generality, 
let $A_1 = 1$. Our basic condition will be that  
$\lambda_1^2 \leq {1}/{3}$, similarly to the first part 
of the proof of Theorem \ref{1}. Note that if $\lambda_1^2 \leq {1}/{3}$ then necessarily  $n \geq 3$.

As easy to check, for $Z \sim N(0,1)$ and $a \in \R$, 
$$
\E\,e^{it\,(Z-a)^2} = \frac{1}{\sqrt{1 - 2it}}\,
\exp\Big\{a^2\,\frac{it}{1 - it}\Big\}, \qquad t \in \R,
$$
so that
$$
\Big|\E\,e^{it\,(Z-a)^2}\Big| = \frac{1}{(1 + 4t^2)^{1/4}}\,
\exp\Big\{-2a^2\,\frac{t^2}{1 + 4t^2}\Big\}.
$$
Hence, the characteristic function $f(t)$ of $W$ satisfies
$$
- \log |f(t)| = 
\frac{1}{4} \sum_{k=1}^n \log(1 + 4\lambda_k^2 t^2)
 + 2 \sum_{k=1}^n 
a_k^2\,\frac{\lambda_k^2 t^2}{1 + 4\lambda_k^2 t^2}.
$$

Since $\lambda_1^2 \leq \frac{1}{3}$, by the monotonicity, all 
$\lambda_k^2 \leq \frac{1}{3}$ as well. But, as we have already 
observed, under the conditons 
$$0 \leq b_k \leq \frac{1}{3}, \ \ \  
b_1 + \dots + b_k = 1,$$ 
and for any fixed value $t \in \R$, the function
$$
\psi(b_1,\dots,b_n) = \sum_{k=1}^n \log(1 + 4b_k t^2)
$$
is minimized for the vector with coordinates 
$$b_1 = b_2 = b_3 = \frac{1}{3} \ \ \ {\rm and}  \ \ \ b_k = 0 \ \ \ 
{\rm for} \ \ \ k>3.$$ 
Hence,
$$
\psi(b_1,\dots,b_n) \geq 3\, \log(1 + 4t^2/3) \geq 
3\, \log(1 + t^2).
$$
Therefore, one may  conclude that
\be
\label{char2}
|f(t)| \leq \frac{1}{(1 + t^2)^{3/4}}\,
\exp\Big\{-2 \sum_{k=1}^n 
a_k^2\,\frac{\lambda_k^2 t^2}{1 + 4\lambda_k^2 t^2}\Big\}.
\en

It is time to involve the inversion formula which yields the upper bound
\be
\label{char3}
M(W) \leq \frac{1}{\pi} \int_0^\infty |f(t)|\,dt.
\en
In the interval 
$$0 < t < T = \frac{1}{2\lambda_1},$$ 
we have 
$\lambda_k^2 t^2 \leq {1}/{4}$ for all $k$, and the bound (8)
is simplified to
$$
|f(t)| \leq \frac{1}{(1 + t^2)^{3/4}}\,e^{-B_1 t^2}.
$$
This gives
$$
\int_0^T |f(t)|\,dt \leq I(B_1) \equiv
\int_0^\infty \frac{1}{(1 + t^2)^{3/4}}\,e^{-B_1 t^2}\,dt.
$$
If $B_1 \leq 1$, 
$$
I(B_1) \leq \int_0^\infty \frac{1}{(1 + t^2)^{3/4}}\,dt < 3,
$$
while for $B_1 \geq 1$,
$$
I(B_1) \leq \int_0^\infty e^{-B_1 t^2}\,dt =
\frac{\sqrt{\pi}}{2\sqrt{B_1}} < \frac{1}{\sqrt{B_1}}.
$$
The two estimates can be united by 
$$I(B_1) \leq \frac{ 3\sqrt{2}}{\sqrt{1+B_1}}.$$
%with $c = 3\sqrt{2}$. 

To perform the integration over the half-axis $t \geq T$,
a different argument is needed.
Put $p_k = a_k^2 \lambda_k^2/B_1$, so that $p_k \geq 0$ and
$p_1 + \dots + p_k = 1$. By Jensen's inequality applied to the
convex function $V(x) = {1}/{(1 + x)}$ for $x \geq 0$ with
points $x_k = 4\lambda_k^2 t^2$, we have
\bee
\sum_{k=1}^n 
a_k^2\,\frac{\lambda_k^2 t^2}{1 + 4\lambda_k^2 t^2}
 & = &
B_1 t^2 \sum_{k=1}^n p_k V(x_k) \\
 & \geq &
B_1 t^2\,V(p_1 x_1 + \dots p_n x_n) \\
 & = &
\frac{B_1 t^2}{1 + \frac{4t^2}{B_1} 
\sum_{k=1}^n a_k^2 \lambda_k^4} \ \geq \
\frac{B_1 t^2}{1 + \frac{4t^2}{3B_1} \sum_{k=1}^n a_k^2 \lambda_k^2} \, = \,
\frac{B_1 t^2}{1 + \frac{4}{3}\,t^2},
\ene
where we used the property $\lambda_k^2 \leq {1}/{3}$.
Moreover, since $$t^2 \geq \frac{1}{(2\lambda_1)^2} \geq \frac{3}{4},$$
necessarily
$$
\frac{t^2}{1 + \frac{4}{3}\,t^2} \geq \frac{3}{8}.
$$
Hence, from (\ref{char2}) we get
$$
|f(t)| \leq \frac{1}{(1 + t^2)^{3/4}}\,e^{-3B_1/4}, \quad t \geq T,
$$
and
$$
\int_T^\infty |f(t)|\,dt \leq e^{-3B_1/4}
\int_{\sqrt{3}/2}^\infty \frac{1}{(1 + t^2)^{3/4}}\,dt < 1.68\, e^{-3B_1/4} < 
\frac{1.85}{\sqrt{1+B_1}}.
$$

Combining the two estimates together for different regions of integration with 
${(3\sqrt{2} + 1.85)}/{\pi} < 1.94$, the bound (\ref{char3}) leads to
$$M(W) < \frac{2}{\sqrt{A_1+B_1}}.$$ 
Thus, this inequality, together with Lemma \ref{lowg}, completes the proof of the theorem.
 
%
%\vskip5mm
%{\bf Proposition 4.} If  $\lambda_1^2 \leq \frac{1}{3} A_1$, then for 
%the maximum of the random variable $W$ defined in (6), 
%there is a two-sided bound
%\be
%\frac{1}{4\sqrt{3}}\,\frac{1}{\sqrt{A_1 + B_1}} \leq
%M(W) \leq \frac{2}{\sqrt{A_1 + B_1}}.
%\en
%Moreover, the left inequality holds true without any assumption.

\section{Acknowledgments}

The research was done within the framework of the Moscow Center for Fundamental and Applied Mathematics, Lomonosov Moscow State
University, and HSE University Basic Research Programs.  Theorem 1 was proved under   support of the RSF grant No. 18-11-00132.. Research of S. Bobkov was supported by the NSF grant
DMS-1855575.

\end{document}